\documentclass [11pt]{amsart}
\usepackage{epic,eepic,latexsym, amssymb, amscd, amsfonts}

\input xy
\xyoption {all}

\def \p {{\mathbb P}}
\def \zh {{{\mathbb Z}/{h\mathbb Z}}}
\def \aut {{\mathsf {Aut} (\mathsf H_h, \mu_{2h})}}
\newtheorem{theorem}{Theorem}
\newtheorem {lemma}{Lemma}

\newtheorem {claim}{Claim} 

\theoremstyle{definition}
\newtheorem{remark}{Remark}

\begin{document}

\title[A note on the Verlinde bundles on elliptic curves]{A note on the Verlinde bundles on elliptic curves}
\author {Dragos Oprea}
\address {Department of Mathematics}
\address {Stanford University}
\email {oprea@math.stanford.edu}
\date{} 

\begin {abstract}We study the splitting properties of the Verlinde bundles over elliptic curves. Our methods rely on the explicit description of the moduli space of semistable vector bundles on elliptic curves, and on the analysis of the  symmetric powers of the Schr\"{o}dinger representation of the Theta group. 
\end{abstract} \maketitle

\section{Introduction} Recently, Popa defined and studied a class of vector bundles on the Jacobians of curves, which he termed the Verlinde bundles \cite {Po1}. The fibers of these vector bundles are the spaces of nonabelian theta functions on the moduli spaces of bundles with fixed determinant over the curve, as the determinant varies in the Jacobian. Popa investigated the splitting properties of these bundles under certain \'{e}tale pullbacks. He further used these results to prove the Strange Duality conjecture at level $1$, and to study the basepointfreeness of the pluri-Theta series. 

In this note, we will study the Verlinde vector bundles in genus $1$. We hope that the results of this work could be useful for the understanding of the higher genus case. In fact, it may be possible to extend our methods to work out a few other low rank/low genus examples. 

To set the stage, consider a smooth complex projective curve $X$ of genus $g\geq 1$, and write $U_{X}(r, r(g-1))$ for the moduli space of rank $r$, degree $r(g-1)$ semistable bundles on $X$.  This moduli space comes equipped with a canonical Theta divisor supported on the locus \begin{equation}\label{thetan}\Theta_{r}=\left\{V \in U_{X}(r, r(g-1)), \text { such that } h^{0}(V)=h^{1}(V)\neq 0\right\}.\end{equation} Following Popa \cite {Po1}, we define the level $k$ Verlinde bundles on the Jacobian as the pushforwards \begin{equation}\label{verl}\mathsf E_{r, k}=\text {det}_{\star}\left(\Theta_{r}^{k}\right)\end{equation} under the determinant morphism $$\text {det}:U_{X}(r, r(g-1))\to \text{Jac}^{r(g-1)}(X)\cong \text{Jac}(X).$$ 

Among the results Popa proved, we mention:
\begin {itemize}
\item [(i)] the pullback of $\mathsf E_{r,k}$ under the multiplication morphism $$r:\text{Jac}(X)\to \text{Jac}(X)$$ splits as a sum of line bundles; 
\item [(ii)] $\mathsf E_{r, k}$ is globally generated iff $k\geq r+1$, and is normally generated iff $k\geq 2r+1$;
\item [(iii)] $\mathsf E_{r, k}$ is ample, polystable with respect to any polarization on the Jacobian, and satisfies $IT_{0}.$
\end {itemize} In addition, it is known that the Verlinde bundles enjoy the following level-rank symmetry:
\begin{itemize}
\item [(iv)]  there is an isomorphism $$\mathsf {SD}:\mathsf E_{r,k}^{\vee}\cong \widehat {\mathsf E_{k,r}}.$$  The hat decorating the bundle on the right hand side denotes the Fourier-Mukai transform with kernel the normalized Poincar\'{e} bundle on the Jacobian. 
\end {itemize}
The morphism $\text{(iv)}$, sometimes termed {\it ``Strange Duality,''} was constructed in this form by Popa. Proofs that $\mathsf {SD}$ is an isomorphism can be found in \cite {mo} \cite {bel}. The case of elliptic curves, which will be relevant for us, is simpler; a discussion is contained in \cite {dt}.
\vskip.1in

To explain the results of this note, assume from now on that $X$ is a smooth {\it complex} projective curve of genus $1$. For reasons which will become clear only later, let us temporarily write $h$ for the rank of the bundles making up the moduli space. We will first show:

\begin{theorem}\label{prop1} Let $k$, $h$ and $q$ be positive integers. The Verlinde bundle $\mathsf E_{h, k}$ splits as a sum of line bundles iff the level $k$ is divisible by the rank $h$. When $k=h(q-1)$, we have \begin{equation}\label{spl}\mathsf E_{h,h(q-1)}\cong\Theta^{q-1}\otimes \left(\bigoplus \mathsf L_{\xi}^{\oplus \mathsf m_{\xi}}\right).\end{equation} Here, $\Theta$ is the canonical Theta bundle on the Jacobian, and the $\mathsf L_{\xi}$'s are the $h$-torsion line bundles. Each line bundle $\mathsf L_{\xi}$ of order $\omega$ occurs with multiplicity \begin{equation}\label{f1}\mathsf m_\xi=\sum_{\delta |h} \frac{1}{q\delta^2}\binom{q \delta}{\delta} \left\{\frac{h/\omega}{h/\delta}\right\},\end{equation} provided that either $h$ or $q$ is odd. If both $h$ and $q$ are even, then \begin{equation}\label{f2}\mathsf m_\xi=\sum_{\delta|h} \frac{(-1)^{\delta}}{q \delta^2}\binom{q \delta}{\delta} \left\{\frac{h/\omega}{h/\delta}\right\}.\end{equation} 
\end {theorem}

The symbol $\{\}$ appearing in the above statement is defined as follows. For any integer $h\geq 2$, we decompose $$h=p_1^{a_1}\ldots p_n^{a_n}$$ into powers of primes. We set \begin{equation}\label{symb}\left\{\frac{\lambda}{h}\right\}=\begin {cases} 0 & \text {if } p_1^{a_1-1}\ldots p_n^{a_n-1} \text {does not divide } \lambda,\\ {\prod}_{i=1}^{n} \left(\epsilon_{i}-\frac{1}{p_i^2}\right) & \text {otherwise }.\end {cases}\end{equation} Here, $$\epsilon_{i}=\begin {cases} 1 & \text{ if } p_{i}^{a_{i}} |\lambda,\\ 0 & \text {otherwise.} \end {cases}$$ 
If $h=1$, the symbol is always defined to be $1$. 
\vskip.1in

Note that it was expected that the splitting of the Verlinde bundles should involve only $h$-torsion line bundles. In fact, Popa proved the isomorphism $$h^{\star}\mathsf E_{h, h(q-1)}\cong h^{\star}\Theta^{N},$$ where $N=\frac{1}{q}\binom{hq}{h}.$ However, the multiplicities $\mathsf m_{\xi}$ of the nontrivial bundles $\mathsf L_{\xi}$ were incorrectly claimed to be $0$ in Propositon $2.7$ of \cite {Po1}. This led to an erroneous statement in Proposition $5.3$. Our note corrects this oversight. 

As an example, when $h$ is an odd prime, all nontrivial $h$-torsion line bundles appear in the decomposition \eqref {spl} with the same (nonzero) multiplicity. This follows for instance by the arguments of \cite {beauville}, upon analyzing the action of a symplectic group on the $h$-torsion points. This is consistent with the Theorem above, which specializes to $$\mathsf E_{h,h(q-1)}=\Theta^{q-1}\otimes \left(\bigoplus_{\xi\neq 0}\mathsf L_\xi^{\oplus n} \oplus \mathcal O^{\oplus m}\right).$$ Here,  $$n=\frac{1}{h^2}\left(\frac{1}{q}\binom{q h}{h}-1\right), \text { and } m=\frac{1}{h^2}\left(\frac{1}{q}\binom{q h}{h}-1\right)+1.$$ Our proof will show that $$m=\dim\left(\text{Sym}^{h(q-1)}\mathsf S_h\right)^{\mathsf H_h},$$ with $\mathsf S_{h}$ being the Schr\"{o}dinger representation of the Heisenberg group $\mathsf H_{h}$. If $h$ is not prime, the ensuing formulas for multiplicities are more complicated, and their integrality is not immediately clear.  
\vskip.1in

Theorem \ref {prop1} is stated for the moduli spaces of bundles of degree zero. The case of arbitrary rank and degree, and of arbitrary Theta divisors will be the subject of Theorem \ref{prop2} in Section \ref{three}.\vskip.07in

The case when the level is not divisible by the rank is slightly more involved, and requires additional ideas. We will consider this most general situation separately, in Section \ref {four}. To explain the final result, let us first change the notation, writing $hr$ for the rank of the bundles making up the moduli space, and letting $hk$ be the level. If $\gcd(r, k)=1$, then, for any $h$-torsion line bunde $\xi$, there is a unique stable bundle $\mathsf W_{r, k, \xi}$ on the Jacobian, having rank $r$ and determinant $\Theta^{k}\otimes \xi$. We will show: 

\begin {theorem}\label{prop3}
Assume that $\gcd (r, k)=1$. The Verlinde bundle of level $hk$ splits as \begin{equation}\label{dec}\mathsf E_{hr,hk}\cong \bigoplus_{\xi}\mathsf W_{r, k, \xi}^{\oplus{\mathsf m}_{\xi}}.\end{equation} For each $h$-torsion line bundle $\xi$ on the Jacobian, having order $\omega$, the multiplicity of the bundle $\mathsf W_{r, k, \xi}$ in the above decomposition equals \begin{equation}\label{mul1}\mathsf m_{\xi}=\sum_{\delta|h}\frac{(-1)^{(h+1)kr\delta}}{(r+k)\delta^{2}} \binom{(r+k)\delta}{r\delta} \left\{\frac{h/\omega}{h/\delta}\right\}.\end{equation} \end {theorem} 

The methods of this work make use of the characteristic zero hypothesis. In positive characteristic, it is likely that the answer is different, and that it depends on the Hasse invariant of the curve. Also, one may justifiably wonder about the higher genus case. This may require a different argument, possibly involving the spaces of conformal blocks. 

\vskip.1in

{\bf Acknowledgements.} We would like to thank Alina Marian for conversations related to this topic and many useful suggestions, and Mihnea Popa for e-mail correspondence. The author was partially supported by the NSF grant DMS-0701114 during the preparation of this work. 

\section{The proof of Theorem \ref{prop1}}\label{two}

The title of this section is self-explanatory. The proof of Theorem \ref {prop1} to be given below relies on two essential ingredients: 
\begin{itemize}
\item [(i)] first, the {\it geometric} input is provided by the explicit description of the moduli space of bundles over elliptic curves, as found in \cite {atiyah}\cite {T};
\item [(ii)] secondly, it will be crucial to understand the symmetric powers of the Schr\"{o}dinger representation of the Heisenberg group.  An {\it algebraic} computation will determine their characters, which are related to the decomposition \eqref{spl} . \end {itemize} 

We will discuss these two items at some length in the next sections, attempting to keep the exposition reasonably self-contained. Our arguments are quite elementary, so it is plausible that some of the results below may already exist in the literature; we tried to provide references, whenever possible. 

\subsection {Geometry}\label{geom} Fix an elliptic curve $(X, o)$. Throughout the paper we identify 
$X\cong \text{Jac}^{0}(X)$ in the usual way, $$p\to \mathcal O_{X}(p-o).$$ In \cite{atiyah}\cite{T}, Atiyah and Tu showed that the moduli space $U_{X}(h, 0)$ of rank $h$, degree $0$ semistable vector bundles on $X$ is isomorphic to the symmetric product $$U_X(h,0)\cong\text{Sym}^h X.$$ Up to $S$-equivalence, the isomorphism can be realized explicitly as \begin{equation}\label{izom}\text{Sym}^h X\ni (p_1,\ldots, p_h)\to \mathcal O_{X}(p_1-o)\oplus \ldots 
\oplus \mathcal O_X(p_h-o)\in U_X(h,0).\end{equation} Under these identifications, the morphism taking bundles to their determinants $$\det:U_X(h,0)\to \text {Jac}(X)$$ is the Abel-Jacobi map, which in this case becomes the addition $$a: \text{Sym}^h X\to X, \, (p_{1}, \ldots, p_{h})\mapsto p_{1}+\ldots+p_{h}.$$ 

Note that the fiber of the morphism $a$ over the point $p\in X$ is the linear series
$$|[p]+(h-1)[o]|=|[p]-[o]+h[o]|.$$ In fact, as an Abel-Jacobi map, the morphism $a$ has the structure of a
projective bundle $\mathbb P(\mathsf V_h)\to X$, where $\mathsf V_h$ is a
rank $h$ vector bundle on $X$.  To describe $\mathsf V_h$, we let $\mathcal P$ be the Poincar\'{e} bundle over
$X\times X$, normalized in the usual way $$\mathcal P=\mathcal O_{X\times
X}(\Delta- \{o\}\times X-X\times \{o\}),$$ with $\Delta\hookrightarrow X \times X$ being the diagonal. Then, using the Fourier-Mukai transform with
kernel $\mathcal P $, denoted $${\bf R} \mathcal S: {\bf  D}(X)\to {\bf D}(X),$$ we have \begin{equation}\label{fmt}\mathsf V_h={\bf R} \mathcal S(\mathcal O_{X}(h [o])).\end{equation} Note that $\mathsf V_h$ has rank $h$, determinant $-[o]$, and, as the Fourier-Mukai
transform of a simple bundle, is simple. In fact, by Atiyah's classic study \cite {atiyah}, there is a unique such bundle on $X$, defined inductively as the (unique) nontrivial extension \begin{equation}\label{fmtt}0\to \mathsf V_{h-1}\to
\mathsf V_{h}\to \mathcal O_{X}\to 0\end{equation} with $\mathsf V_1=\mathcal O_{X}(-[o]).$ Alternatively, this exact sequence is obtained as the Fourier-Mukai transform of $$0\to \mathcal O_{X}((h-1)[o])\to \mathcal O_{X}(h[o])\to \mathcal O_{\{o\}}\to 0.$$

Note that the line bundle \eqref{izom} has a section precisely when $p_{i}=o$ for some $1\leq i\leq h$. It follows from $\eqref{thetan}$ that the canonical theta divisor $\Theta_{h}$ on $U_{X} (h,0)$ is the image of the symmetric 
sum \begin{equation}\label{thd}[o]+\text{Sym}^{h-1}X\hookrightarrow \text{Sym}^h X.\end{equation} Thus, the Theta line bundle $\Theta_{h}$ agrees, at least
fiberwise, with $\mathcal O_{\mathbb P(\mathsf V_h)}(1)$.  In fact, one can show the isomorphism $$\Theta_{h}\cong\mathcal O_{\mathbb P(\mathsf V_{h})}(1).$$ Moreover, the canonical section vanishing along the Theta divisor \eqref{thd} is the composition $$\mathcal O_{\p(\mathsf V_{h})}\to \mathcal O_{\p(\mathsf V_{h})}(1)\otimes a^{\star}\mathsf V_{h}\to \mathcal O_{\p(\mathsf V_{h})}(1),$$ with the second arrow given by \eqref{fmtt}. These observations allow us to compute the level $k$ Verlinde bundle \begin{equation}\label{ehk}\mathsf E_{h,k}=a_{\star}\left(\Theta_{h}^{k}\right)=a_{\star}
\left(\mathcal O_{\mathbb P(\mathsf V_{h})}(k)\right)=\text {Sym}^k \mathsf V_h^{\vee}.\end{equation}

For convenience, we will write $\mathsf W_h=\mathsf V_h^{\vee}$ for the unique stable bundle on $X$ of rank $h$ and determinant $\mathcal O_{X}([o])$. More generally, if $\gcd (h, d)=1$, we let $\mathsf W_{h, d}$ be the unique stable bundle of rank $h$ and determinant $\mathcal O_{X}(d[o])$. The bundles $\mathsf W_{h,d}$ can be constructed inductively as successive extensions \cite {polishchuck}. Indeed, consider two consecutive terms $0\leq \frac{d_{1}}{h_{1}}<\frac{d_{2}}{h_{2}}<1$ in the Farey sequence, {\it i.e.} assume that $$h_{1}d_{2}-h_{2}d_{1}=1.$$ Set $h=h_{1}+h_{2}, d=d_{1}+d_{2}$. Then $\mathsf W_{h,d}$ is the unique nontrivial extension $$0\to \mathsf W_{h_{1}, d_{1}}\to \mathsf W_{h, d}\to \mathsf W_{h_{2}, d_{2}}\to 0.$$

\vskip.1in

With these preliminaries out of the way, we proceed to investigate the splitting behavior of the Verlinde bundles $\mathsf E_{h,k}$. Our analysis relies on the multiplicative structure of the Atiyah bundles \cite {atiyah}, which may not be immediately obvious.

\begin{lemma}\label{l1} The Verlinde bundle $\mathsf E_{h,k}$ splits as a sum of
line bundles if and only if $h$ divides $k$. \end {lemma}

\proof This result will be reproved later in the paper. A more direct argument is given below. First, observe that $\mathsf E_{h,k}$ is a direct summand of $\mathsf W_h^{\otimes k}$. It suffices to show that these tensor powers split as sums of lines bundles iff $h$ divides $k$. 
In fact, something more general is true: 
\begin{claim}
Assuming $\gcd (h, d)=1$, the tensor powers $\mathsf W_{h,d}^{\otimes k}$ split as sums of rank $h'$ bundles of the form $\mathsf W_{h', d k'}\otimes M$ where $M$ are various degree $0$ line bundles. Here, we set $$h'=\frac{h}{\gcd(h,k)}, k'= \frac{k}{\gcd(h,k)}.$$
\end {claim}

To prove the {\it Claim}, we first decompose $h=h_1\ldots h_s$ into powers of primes, and pick integers $d_1, \ldots, d_s$ such that $$\frac{d_1}{h_1}+\ldots+\frac{d_s}{h_s}=\frac{d}{h}.$$ Then, $$\mathsf W_{h,d}=\mathsf W_{h_1, d_1}\otimes \ldots \otimes \mathsf W_{h_s,d_s}.$$ This could be argued as follows: both sides have the same (coprime) rank and determinant, and are moreover semistable, in fact stable. Therefore, they should coincide by Atiyah's classification. With this understood, one checks that it is enough to take $h$ to be a power of a prime $p$. 

For the latter case, we will need the following rephrasing of Theorems $13$ and $14$ in \cite {atiyah}. Assume $e_{1}, e_{2}, e$ are integers not divisible by $p$, and that $$\frac{e_{1}}{p^{a_{1}}}+\frac{e_{2}}{p^{a_{2}}}=\frac{e}{p^{a}}.$$ Then, Atiyah showed that for certain degree $0$ line bundles $M$, we have  \begin{equation}\label{at}
\mathsf W_{p^{a_{1}}, e_{1}}\otimes \mathsf W_{p^{a_{2}}, e_{2}}=\bigoplus_{M}\mathsf W_{p^{a}, e}\otimes M.
\end {equation}

Thus, when $h$ is a power of a prime, the {\it Claim} follows from \eqref{at}, by a straightforward induction on $k$. \qed

\begin {remark} Using a sharper version of Atiyah's results, one can prove that when $h$ is odd, the $M$'s appearing in the {\it Claim} above are representatatives for the cosets of $h$-torsion line bundles on $X$ modulo the twisting action of the group of $h'$-torsion line bundles. The same statement should hold true for $h$ even, but Atiyah's results only show that the orders of the $M$'s divide $2h$. In particular, for $h$ odd and $\gcd (h, k)=1$, we immediately conclude that \begin{equation}\label{ehk}\mathsf E_{h,k}\cong\bigoplus_{i=1}^{m}\mathsf W_{h, k},\end{equation} with $m=\frac{1}{h+k}\binom{h+k}{h}.$ Equation \eqref{ehk} is a particular case of Theorem \ref {prop3}.
\end {remark}

\vskip.1in
We will identify the splitting of $\mathsf E_{h,k}=\text{Sym}^{k}\mathsf W_h$ when the level $k$ is divisible by the rank $h$. We set $$q=1+\frac{k}{h}.$$ 

Let $\mathsf X_h$ be the group of $h$-torsion points on the elliptic curve. Let $\mathsf G_h$ be the Theta group of the line bundle $\mathcal O_{X}(h[o]),$ which is a central extension $$1\to \mathbb C^{\star}\to \mathsf G_h\to \mathsf X_h\to 1.$$ The assignment $$\eta\to \eta^{2h}$$ defines an endomorphism of $\mathsf G_{h}$, whose image lies in the center of $\mathsf G_{h}$. Let $\mathsf H_{h}$ be the kernel of this endomorphism. It corresponds to an extension $$1\to \mu_{2h} \to \mathsf H_h\to \mathsf X_h\to 1,$$ where $\mu_{2h}\hookrightarrow \mathbb C^{\star}$ is the group of $2h$-roots of $1$. Finally, let $\mathsf S_h$ denote the $h$-dimensional Schr\"{o}dinger representation of $\mathsf G_h$, {\it i.e.} the unique representation such that the center of $\mathsf G_{h}$ acts by its natural character. 

Picking theta structures, we identify $\mathsf G_{h}$ with the Heisenberg group \begin{equation}\label{ident}\mathsf G_h\cong \mathbb C^{\star} \times \zh\times \zh.\end{equation} The multiplication on the right hand side is defined as $$(\alpha, x, y) (\alpha', x', y')=(\alpha \alpha' \zeta^{y'x}, x+x', y+y').$$ Here, we set $$\zeta=\exp \left(\frac{2\pi i}{h}\right).$$ The Schr\"{o}dinger representation $\mathsf S_h$ is realized on the space of functions $$f:\zh \to \mathbb C.$$ The action of the element $(\alpha, x, y)\in \mathsf G_{h}$ on a function $f$ is given by the new function $$F:\zh \to \mathbb C, \,\, F(a)=\alpha \zeta^{ya}\cdot f(x+a).$$ 
\vskip.1in

We will first compute the pullbacks of the Verlinde bundles under the morphism $h:X\to X$ which multiplies by $h$ on the elliptic curve. Using the description of $\mathsf V_{h}$ as a Fourier-Mukai transform provided by \eqref{fmt}, and Theorem $3.11$ in \cite {mukai}, we obtain $$h^{\star}\mathsf V_{h}\cong \mathcal O_{X}(-h[o])^{\oplus h}.$$ In fact, we claim that $\mathsf G_{}{h}$-equivariantly, we have \cite{polishchuck} \begin{equation}\label{ess}h^{\star} \mathsf V_h\cong \mathsf S_h\otimes \mathcal O_{X}(-h[o]).\end{equation} Indeed, consider the trivial bundle $$h^{\star}\mathsf V_{h}\otimes \mathcal O_{X}(h[o])\cong V\otimes \mathcal O_{X},$$ where $V$ is an $h$-dimensional vector space. Both factors of the tensor product on the left carry a $\mathsf G_{h}$-action covering the translation $\mathsf X_{h}$-action on the base $X$.  Therefore, endowing the structure sheaf appearing on the right with the trivial $\mathsf G_{h}$-action, we obtain an $\mathsf G_{h}$-representation on $V$.  Moreover, note that the center of $\mathsf G_{h}$ acts on $V$ by homotheties. Therefore, $V\cong \mathsf S_{h},$ by the uniqueness of the Schr\"{o}dinger representation. This establishes \eqref{ess}. 

Taking determinants in \eqref{ess}, we obtain \begin{equation}\label{dets}h^{\star}\mathcal O_{X}(-[o])\cong \Lambda^{h}\mathsf S_{h}\otimes \mathcal O_{X}(-h[o])^{h}.\end{equation} This identification is a priori only $\mathsf G_{h}$-equivariant, but, since the center of $\mathsf G_{h}$ acts trivially, the isomorphism is in fact $\mathsf X_{h}$-equivariant. Similarly, dualizing and taking symmetric powers in \eqref{ess}, we obtain an $\mathsf X_h$-equivariant identification \begin{equation}\label{lat1}h^{\star} \text{Sym}^{k} \mathsf W_h\cong\text{Sym}^{k}\mathsf S_h^{\vee}\otimes \mathcal O_{X}(h[o])^{k}\cong \text{Sym}^{k}\mathsf S_{h}^{\vee}\otimes \left(\Lambda^{h}\mathsf S_{h}\right)^{q-1}\otimes h^{\star}\mathcal O_{X}([o])^{q -1}.\end{equation} 

Observe that the action of the central elements $\alpha$ of $\mathsf G_h$ on the Heisenberg module $$\mathsf M_k=\text{Sym}^{k}\mathsf S_{h}^{\vee}\otimes \left(\Lambda^{h}\mathsf S_{h}\right)^{q-1}$$  is trivial, since \begin{equation}\label{scaling}\alpha^{-k}\cdot (\alpha^{h})^{q-1}=1.\end{equation}  Therefore $\mathsf M_{k}$ is  an $\mathsf X_{h}$-module. The $\mathsf X_{h}$-action splits into eigenspaces indexed by the characters $\xi$ of $\mathsf X_{h}$, each appearing with multiplicity $\mathsf m_\xi$: \begin{equation}\label{decomp}\mathsf M_{k}\cong \bigoplus_{\xi}\xi^{\oplus \mathsf m_{\xi}}.\end{equation}

Let us write $\widehat{\mathsf X}_h$ for the group of characters of $\mathsf X_h$. For each $\xi\in \widehat{\mathsf X}_h$, let $\mathsf L_\xi$ denote the corresponding $h$-torsion line bundle on $X$. The pullback $h^{\star}\mathsf L_{\xi}$ is the trivial bundle endowed with the $\mathsf X_{h}$-character $\xi$. Using \eqref{lat1} and \eqref{decomp}, we obtain an $\mathsf X_{h}$-equivariant identification $$h^{\star} \text{Sym}^{k} \mathsf W_h\cong h^{\star}\left(\bigoplus_{\xi}\mathsf L_{\xi}^{\oplus\mathsf m_{\xi}}\right)\otimes h^{\star}\mathcal O_{X}([o])^{q-1}.$$ 

This equivariant isomorphism determines the Verlinde bundle on the left, by general considerations about the Picard group of finite quotients. We can also give a direct argument as follows. Pushing forward the previous equation by $h$, we obtain the $\mathsf X_{h}$-isomorphism \begin{equation}\label{push1}\text{Sym}^{k}\mathsf W_{h}\otimes h_{\star}\mathcal O_{X}\cong \bigoplus_{\xi}\mathsf L_{\xi}^{\oplus \mathsf m_{\xi}} \otimes \mathcal O_{X}([o])^{q-1}\otimes h_{\star}\mathcal O_{X}.\end {equation} Note that $\mathsf X_h$-equivariantly \begin{equation}\label{push}h_\star \mathcal O_{X}\cong\sum_{\xi\in\widehat{\mathsf X}_h} \mathsf L_{\xi}.\end {equation} Comparing \eqref{push1} and \eqref{push}, and singling out the $\mathsf X_{h}$-invariant part, we conclude that \begin{equation}\label{split}\mathsf E_{h,k}\cong\text{Sym}^{k} \mathsf W_h\cong\left(\bigoplus_{\xi\in \widehat{\mathsf X}_h} \mathsf L_{\xi}^{\oplus \mathsf m_{\xi}}\right)\otimes \mathcal O_{X}\left([o]\right)^{q-1}.\end{equation}

\subsection {Algebra} \label{alg} It remains to determine the multiplicities $\mathsf m_\xi$ appearing in \eqref{decomp}. Regarding $\mathsf M_k$ as a representation of the finite group $\mathsf H_h$, it is clear that \begin{equation}\label{sum}\mathsf m_\xi=\frac{1}{|\mathsf H_h|}\sum_{\eta\in \mathsf H_h} \xi(\eta^{-1}) \text{Tr}_{\mathsf M_{k}}(\eta).\end{equation} We will compute this sum explicitly with the aid of the following 

\begin{lemma}\label{chr} Let $\eta\in \mathsf H_h$ be an element whose image under the map $\mathsf H_h\to \mathsf X_h$ has order exactly $h/\delta$ in $\mathsf X_h$. The trace of $\eta$ on $\mathsf M_k$ equals $$\text{Tr} _{\mathsf M_{k}}(\eta)=\frac{1}{q}\binom{q \delta}{\delta},$$ provided that either $h$ or $q$ is odd. 

\end {lemma}

\proof  We pick theta structures, so that $\mathsf G_h$ and $\mathsf S_h$ are given by \eqref{ident}. Consider the basis $f_1, \ldots, f_h$ of $\mathsf S_h$ given by $$f_i(j)=\delta_{i,j}.$$ By definition, the action of $$\eta=(\alpha, x,y)\in \mu_{2h}\times \zh \times \zh$$ is given as \begin{equation}\label{act}\eta \cdot f_i=\alpha \zeta^{y(i-x)}\cdot f_{i-x}.\end{equation} To compute the trace of $\eta$ on $\mathsf M_{k}$, we may assume that $\alpha=1$, since the scaling action of the center of $\mathsf H_{h}$ is trivial, as remarked in \eqref{scaling}. 

We begin by computing the trace $\text{Tr Sym}^{k}\eta$ of the action of $\eta$ on $\text{Sym}^{k}\mathsf S_{h}^{\vee}$. For simplicity, we will first treat the case $x=0$. The eigenvalues of the action of $\eta$ on $\mathsf S_{h}$ are $1, \zeta_{y}, \ldots, \zeta_{y}^{h-1}.$ Here, we set $$\zeta_y=\zeta^{y}.$$ Therefore, \begin{equation}\label{sum1}\text{Tr } \text{Sym}^{k}\eta=\sum_{1\leq i_1\leq \ldots \leq i_k\leq h} \zeta_{y}^{-(i_1+\ldots+i_k)}=\sum_{j_1+\ldots+j_h=k} \zeta_{y}^{-(j_1+2j_2+\ldots+hj_h)}.\end{equation} In the above, $j_r$ denotes the number of $i$'s which equal $r$. Now, we compute the generating series $$\sum_{k}\text{Tr } \text{Sym}^{k}\eta\cdot t^{k}=\frac{1}{1-\zeta_{y}^{-1} t}\cdot \frac{1}{1-\zeta_{y}^{-2}t}\cdots \frac{1}{1-\zeta_{y}^{-h}t}.$$ Write $$h=lm,$$ where $$m=\gcd(h, y) \text { and } \gcd (l,y)=1.$$ Then $\epsilon=\zeta_{y}$ is a primitive root of $1$ of order $l$. Therefore, the product in the denominator above becomes $$(1-\zeta_{y}^{-1} t)\ldots (1-\zeta_{y}^{-h}t)=\zeta_{y}^{-h(h-1)/2}(-1)^h \left((t-1)(t-\epsilon)\ldots (t-\epsilon^{l-1})\right)^{m}$$ $$=(-1)^{h+y(h-1)}(t^{l}-1)^{m}.$$ We can extract the coefficient of $t^{k}$:
\begin{equation}\label{tr}\text{Tr } \text{Sym}^{k}\eta=(-1)^{h+y(h-1)+m+\frac{k}{l}}\binom{-m}{\frac{k}{l}}=(-1)^{\frac{k}{l}}\binom{-m}{\frac{k}{l}}=\frac{1}{q}\binom{q m}{m}.\end{equation} In particular, this computation implies that the sum \eqref{sum1} is $1$ when $m=\gcd (h, y)=1$. Moreover, the argument shows that the sum \eqref{sum1} vanishes if $k$ is not divisible by $l=\frac{h}{\gcd (h,y)}$. 

We will now consider the $\eta$'s in $\mathsf H_h$ for which $x\neq 0$. For these, the computation is notationally more involved. To begin, we write $$x=x's, \text { and } h=h's,$$ where $s=\gcd (h,x)$. Let $u$ be any constant with $$u^{h'}=(-1)^{yx'(h'+1)}.$$ Note in particular that $u^k=1$ for $h$ odd. For $h$ even, we have \begin{equation}\label{sign1}u^{k}=(-1)^{xy(q -1)}.\end{equation}

Now, it is easy to see that the eigenvalues of $\eta$ on $\mathsf S_{h}$ are \begin{equation}\label{egv}\lambda_{i,j}=u \,\zeta_{y}^{i} \,\sigma^{j},\,\, \, 1\leq i\leq s, 1\leq j\leq h',\end {equation} where $$\sigma=\exp\left(\frac{2\pi i}{h'}\right).$$ In fact, we can exhibit an eigenvector for $\lambda=\lambda_{i,j}$, namely $$v_\lambda=\sum_{k=0}^{h'-1}\lambda^{-k}\zeta_y^{ki-\frac{k(k+1)}{2}x}\cdot f_{i-kx}.$$ We order the indices $(i,j)$ lexicographically. The trace $\text{Tr }\text{Sym}^k \eta$ is obtained by summing all products $$\left(\lambda_{1, j^1_1}^{-1}\cdots \lambda_{1, j^1_{\bullet}}^{-1}\right)\left( \lambda_{2, j^2_1}^{-1}\cdots \lambda_{2, j^{2}_{\bullet}}^{-1}\right) \cdots \left(\lambda_{s, j^{s}_1}^{-1}\cdots \lambda_{s, j^{s}_{\bullet}}^{-1}\right),$$ where $$1\leq j^{i}_1\leq j^{i}_2\leq \ldots \leq j^{i}_{\bullet}\leq h'.$$ Let $a_1$ be the number of terms in the product whose first index is $1$; $a_2, \ldots, a_s$ have the similar meaning. We require $a_1+\ldots+a_s=k$. After substituting \eqref{egv} in the product above, we sum over the $j$'s, keeping the $a$'s fixed. We have seen already in the derivation of \eqref{sum1} that the sum $$\sum_{1\leq j^i_1\leq \cdots\leq j^i_{a_i}\leq h'} \sigma^{-(j^i_1+\ldots +j^i_{a_i})}$$ is $0$ if $h'$ does not divide $a_{i}$, and it equals $1$ otherwise. Therefore, writing $a_i=h' a_i'$, we need to evaluate $$\sum_{a'_1+\ldots+a'_s=\frac{k}{h'}} \zeta_y^{-h'a'_1}\zeta_y^{-2h'a'_2}\cdots \zeta_y^{-sh'a'_s}=\sum_{a'_1+\ldots+a'_s=\frac{k}{h'}}\gamma_y^{-(a'_1+\ldots+sa'_s)}.$$ Here, we set $\gamma=\exp\left(\frac{2\pi i}{s}\right)$, so that $\zeta_{y}^{h'}=\gamma_y$. We have already computed sums of this type in \eqref{sum1}. We obtained the answer \begin{equation}\label{ans1}\frac{1}{q}\binom{q \delta}{\delta}\end {equation} for $\delta=\gcd (s ,y)=\gcd (h, x, y).$ This expression gives the trace $\text {Tr Sym}^{k}\eta$ when $h$ is odd. The formula includes the previously considered case $x=0$, for which $\delta=m.$ The sign change \eqref{sign1} is required when $h$ is even. 

Finally, the trace of $\eta$ on $\Lambda^{h}\mathsf S_{h}$ is computed using \eqref{act}: \begin{equation}\label{signs}\eta\cdot f_{1}\wedge\ldots \wedge f_{h}=(-1)^{x(h+1)}\prod_{i=1}^{h}\zeta_{y}^{i-x}\cdot f_{1}\wedge \ldots \wedge f_{h}=(-1)^{(h+1)(x+y)}f_{1}\wedge \ldots \wedge f_{h}.\end{equation} 

This completes the proof when $h$ is odd. When $h$ is even, we take into account the sign corrections of the previous paragraph and  \eqref{sign1}. We append formula \eqref{ans1} by the overall sign $$(-1)^{xy(q-1)}\cdot (-1)^{(x+y)(h+1)(q-1)}=(-1)^{(xy+x+y)(q-1)}.$$ This does not change \eqref{ans1} when $q$ is odd, proving the Lemma. When $h$ and $q$ are both even, we note, for further use, that the overall sign of \eqref{ans1} can be rewritten as \begin{equation}\label{oall}(-1)^{\gcd (h, x, y)}=(-1)^{\delta}.\end{equation}\qed

\vskip.15in

We proceed to calculate the sum \eqref{sum}. We {\it claim} that the multiplicity $\mathsf m_{\xi}$ depends only on the order of the character $\xi\in \widehat{\mathsf X}_h$. To this end, consider the group $\aut$ of automorphisms of $\mathsf H_h$ which restrict to the identity on the center $\mu_{2h}$. As essentially remarked in \cite {beauville},  the characters appearing in the $\mathsf X_{h}$-representation $\mathsf M_{k}$ are exchanged by the action of $\aut$. Beauville's argument is based on the observation that for each $F\in \aut$, the standard $\mathsf H_{h}$-module structure of $\mathsf S_{h}$, $\rho: \mathsf H_{h}\to \text{GL} (\mathsf S_{h}),$ is isomorphic to the twisted module structure $F\circ \rho: \mathsf H_{h}\to \text{GL}(\mathsf S_{h})$. This follows by examining the character of the center of $\mathsf H_{h}$, and by making use of the uniqueness of the Schr\"{o}dinger representation. The same observation applies to the  associated $\mathsf H_{h}$-module $\mathsf M_{k}$. With this understood, our {\it claim} is a consequence of the Lemma below. This result is possibly known, yet for completeness we decided to include the argument. Note that the Lemma is not indispensable for the proofs to follow, yet it allows for some simplification of the formulas.  

\begin {lemma} Under the action of $\aut$, two characters of $\mathsf X_h$ belong to the same orbit if and only if they have the same order in $\widehat {\mathsf X}_h.$ 
\end {lemma} 

\proof Fix two characters $\chi_1, \chi_2$ of $\mathsf X_{h}$: $$\chi_{i}:\mathsf X_h \to \mathbb C^{\star}, (x,y)\to \zeta^{a_ix+b_iy}, \, 1\leq	 i\leq 2.$$ The condition on the orders of $\chi_{1}$ and $\chi_{2}$ translates into $$\gcd(h, a_1, b_1)=\gcd(h, a_2, b_2):=\tau.$$ This implies that we can solve the equations below, with the Greek letters as the unknows: \begin{equation}\label{sys1}a_1 \lambda+b_1 \mu=a_2 \mod h, \, \, a_1 \nu+b_1 \gamma = b_2 \mod h.\end{equation} We claim that we may further achieve \begin{equation}\label{sys2}\lambda \gamma - \mu \nu=1\mod h.\end{equation} This can be seen for instance as follows. By the Chinese Remainder Theorem, we may take $h$ to be a power of a prime. In this case, assume first that $\tau=1$. Starting with any solution of \eqref{sys1}, define a new quadruple $$\lambda'=\lambda+b_{1}x,\, \mu'=\mu-a_{1}x,\,\,  \nu' = \nu+b_{1}y, \,  \gamma'=\gamma-a_{1}y.$$ The assumption $\tau=1$ implies that we can find a pair $(x,y)$ such that \eqref{sys2} holds: $$\lambda'\gamma'-\mu'\nu'=(\lambda \gamma-\mu\nu)+ b_{2}x-a_{2}y\equiv 1\mod h$$ For arbitrary $\tau$, after dividing by $\tau$, and using the case we already proved, we may assume that \eqref{sys1} is satisfied $\mod h$, and that \eqref{sys2} holds true $\mod h/\tau$. We lift the solution using Hensel's lemma, ensuring that \eqref{sys2} is also satisfied $\mod h$. 

Finally, define $F:\mathsf H_h \to \mathsf H_h$ by $$F(\alpha, x, y)=(\alpha \zeta^{\frac{1}{2}({\lambda \mu}x^2+\nu\gamma y^2+2\mu \nu xy)} , \lambda x+ \nu y, \mu x+ \gamma y).$$ Equation \eqref{sys2} is used to prove that $F$ is an automorphism of $\mathsf H_{h}$, while equation \eqref{sys1} shows that $F$ sends $\chi_1$ to $\chi_2$. \qed

\vskip.1in
Henceforth, for the computation of \eqref{sum}, we will take $\xi$ to be the character $$\xi=\xi_{\lambda}:\zh\times \zh \ni (x,y)\mapsto \zeta_\lambda^{x+y}=\zeta^{\lambda(x+y)}\in \mathbb C^{\star}.$$ Here, we assume that $\lambda$ divides $h$, so that the character $\xi$ has order $${\omega}=\frac{h}{\lambda}.$$ Assume that either $h$ or $q$ is odd. Using Lemma \ref{chr}, we rewrite \eqref{sum} as \begin{equation}\label{sum2}\mathsf m_{\xi}=\frac{1}{h^2}\sum_{\delta|h} \frac{1}{q}\binom{q \delta}{\delta} \left(\sum_{\gcd (h, x, y)=\delta} \xi_{\lambda}(x,y)\right).\end {equation} If both $h$ and $q$ are even, each term in \eqref{sum2} is multiplied by the sign $(-1)^{\delta}$, as it follows from \eqref{oall}. In this case, \begin{equation}\label{sum3}\mathsf m_{\xi}=\frac{1}{h^2}\sum_{\delta|h} \frac{(-1)^{\delta}}{q}\binom{q \delta}{\delta} \left(\sum_{\gcd (h, x, y)=\delta} \xi_{\lambda}(x,y)\right).\end {equation}

We will evaluate formulas \eqref{sum2} and \eqref{sum3} in terms of the character \eqref{symb} defined in the introduction.  

\begin {lemma}\label{nl} We have $$\sum_{\gcd (h, x, y)=\delta} \xi_{\lambda}(x, y)=\frac{h^{2}}{\delta^{2}}\left\{\frac{h/\omega}{h/\delta}\right\}.$$
\end {lemma}

\proof Replacing $h$, $x$ and $y$ by $h/\delta$, $x/\delta$ and $y/\delta$ respectively, we may assume $\delta=1$. To solve this case, let us set \begin{equation}\label{n}\mathsf N_{\lambda}(h)=\sum_{\gcd (h, x, y)=1}\xi_{\lambda}(x, y)=\sum_{\gcd (h, x, y)=1}\zeta_\lambda^{x+y}.\end {equation} It suffices to show that \begin{equation}\label{cl}\mathsf N_{\lambda}(h)=h^2 \left\{\frac{\lambda}{h}\right\}.\end{equation}

This is immediate when $h=p^a$ is a power of a prime. In this case, if $p^a|\lambda$, the left hand side of \eqref{cl} counts the pairs $1\leq x,y\leq p^a$ such that $\gcd (p^a, x, y)=1$. Their number is $p^{2a-2}(p^2-1),$ which equals the right hand side. Otherwise, since the distinct roots of unity add up to $0$, we have $$\sum_{(x,y, p^a)=1} \zeta_{\lambda}^{x+y}=-\sum_{p |(x,y)} \zeta_{\lambda}^{x+y}.$$ If $p^{a-1}|\lambda$, then all terms in the last sum are equal to $1$, hence giving the answer $-p^{2a-2}$. Finally, if $p^{a-1}$ does not divide $\lambda$, then replacing $\zeta_{\lambda}$ by $\zeta_{p\lambda}$, we sum all distinct roots of unity of order $p^{a-1}/\gcd(p^{a-1}, \lambda)$, each appearing with equal multiplicity. This gives the answer $0$. 

The general case follows by induction on the number of prime factors of $h$, once we establish the multiplicativity  in $h$ of the function $\mathsf N_{\lambda}(h)$. Let $h=h_1 h_2$ with $\gcd(h_1, h_2)=1$. Chose integers $u,v$ such that $$h_1 u + h_2 v=1.$$ By the Chinese Remainder Theorem, the pairs $(x, y)\mod h$ are in one-to-one correspondence with pairs $(x_1, y_1)\mod h_1$, $(x_2, y_2)\mod h_2$ such that $$x\equiv x_1 \mod h_1,\,\, x\equiv x_2 \mod h_2,$$ $$y\equiv y_1 \mod h_1,\,\, y\equiv y_2 \mod h_2.$$ Explicitly, we have $$x= h_1 u x_2 + h_2 v x_1 \mod h, \,\, y=h_1 u y_2+h_2 v y_1 \mod h.$$ The condition $\gcd (h, x, y)=1$ is equivalent to $$\gcd (h_{1}, x_1, y_1)=1, \, \gcd(h_{2}, x_2, y_2)=1.$$ We compute \begin{eqnarray*} \mathsf N_{\lambda}(h)&=&\sum_{\gcd (h, x, y)=1} \zeta_{\lambda}^{x+y}=\sum_{\gcd (h_1, x_1, y_1)=1, \gcd (h_2, x_2, y_2)=1} \zeta_{\lambda}^{h_2 v(x_1+y_1)}\cdot\zeta_{\lambda}^{h_1 u(x_2+y_2)}\\ &=&\mathsf N_{\lambda v} (h_1)\mathsf N_{\lambda u}(h_2)=h_1^2\left\{\frac{\lambda v}{h_1}\right\}\cdot h_{2}^{2}\left\{\frac{\lambda u}{h_2}\right\}=h^2\left\{\frac{\lambda}{h_1}\right\}\left\{\frac{\lambda}{h_2}\right\}=h^2\left\{\frac{\lambda}{h}\right\}.\nonumber\end {eqnarray*} In the last line, we used the fact that the factors $u$ and $v$ do not change the symbol $\left\{\right\}$ since these numbers are prime to $h_2$ and $h_1$ respectively. 
\qed
\vskip.1in

Putting together \eqref{split}, \eqref{sum2}, \eqref{sum3} and Lemma \ref{nl}, we complete the proof of Theorem \ref{prop1}.\vskip.1in

\section {Arbitrary numerics}

\subsection {Arbitrary rank and degree}\label{three} We will now discuss a variant of Theorem \ref{prop1}, which covers the case of arbitrary rank and degree. Let $r, d$ be two integers with $$h=\gcd(r, d).$$ Write $$r=h r', d=h d', \text { where } \gcd (r', d')=1.$$ We will consider Theta divisors  on the moduli space $U_{X}(r, d)$. Their definition requires the choice of a twisting vector bundle $N$ of complementary slope $$\mu(N)=-\frac{d}{r}.$$ We set \begin{equation}\label{th2}\Theta_{r,N}=\{V\in U_{X} (r, d), \text { such that } h^{0}(V\otimes N)=h^{1}(V\otimes N)\neq 0\}.\end{equation} To avoid confusion, even though it may be notationally cumbersome, we decorate the Theta's by the twisting bundles $N$, and by the rank of the bundles in the moduli space. 

It is convenient to assume that $N$ has the minimal possible rank $r'$. The level $k$ Verlinde bundle $$\mathsf E^{N}_{r, k}={\det}_{\star}\left(\Theta_{r,N}^{k}\right)$$ is obtained by pushing forward the pluri-Theta bundle $\Theta_{N}^{k}$ on $U_{X}(r, d)$ via the morphism $$\det:U_{X}(r, d)\to \text{Jac}^{d}(X).$$ 

As before, we have an isomorphism \begin{equation}\label{arb} U_{X} (r, d)\cong \text{Sym}^{h} X.\end {equation} Set-theoretically, this isomorphism is essentially defined twisting \eqref{izom} by the unique idecomposable vector bundle $\mathsf W_{r',d'}$ of rank $r'$ and determinant $d'[o]$ on $X$. More precisely, if $(p_{1}, \ldots, p_{h})$ are $h$ points of $X$, pick $(q_{1}, \ldots, q_{h})$ such that $$r' \cdot q_{i}=p_{i}, \, 1\leq i\leq h.$$ Then, the isomorphism \eqref{arb} is given by   
\begin{equation}\label{izomex}\text{Sym}^{h} X \ni (p_{1}, \ldots, p_{h})\mapsto \mathsf W_{r', d'}\otimes \mathcal O_{X}(q_{1}-o) \oplus \ldots \oplus \mathsf W_{r', d'}\otimes \mathcal O_{X}(q_{h}-o)\in U_{X} (r, d).\end{equation} Note that the answer on the right hand side of \eqref{izomex} is independent of the choice of $q_{i}$. Indeed, any two $q_{i}$'s must differ by an $r'$-torsion point $\chi$. However, by Atiyah's classification, \begin{equation}
\label{inv}\mathsf W_{r', d'}\otimes \mathsf L_{\chi}\cong \mathsf W_{r', d'},\end{equation} as both bundles are simple, of the same rank and determinant. It was observed in \cite {T}, and it is clear from \eqref{izomex}, that the determinant $$\det: U_{X} (r, d)\to \text{Jac}^{d}(X)$$ becomes the addition morphism $$a: \text{Sym}^{h}X\to X, \,(p_{1}, \ldots, p_{h})\to p_{1}+\ldots+p_{h}.$$ Here, we used the identification $$X\cong \text {Jac}(X)\cong \text {Jac}^{d}(X),$$ with the second arrow given by twisting degree zero line bundles by $\mathcal O_{X}(d[o]).$ Via this identification, the divisor $\Theta_{1, \mathcal O(-d[o])}$ on $\text{Jac}^{d}(X)$ corresponds to the canonical Theta on $\text{Jac} (X)$. 

Finally, we can easily identify the Theta divisors on $U_{X}(r, d)$. There is a natural choice for the twisting bundle $N$, namely the Atiyah bundle $N_{0}=\mathsf W_{r', -d'}.$ It was shown in \cite {atiyah}, and it follows from equation \eqref{inv}, that the tensor product \begin{equation}\label{inv2}\mathsf W_{r', -d'}\otimes \mathsf W_{r', d'}=\bigoplus_{\chi}\mathsf L_{\chi}\end{equation} splits as the direct sum of all $r'$-torsion line bundles $\mathsf L_{\chi}.$ As a consequence of \eqref{th2}, \eqref{izomex}, \eqref{inv2}, we see that for the bundles $V$ in the Theta divisor, we have $q_{i}=\chi$ for some $r'$-torsion point $\chi$, and some $1\leq i\leq h$. Thus, $\Theta_{r, N_{0}}$ is the image of the symmetric sum $$[o]+\text{Sym}^{h-1}X\hookrightarrow \text{Sym}^{h}X.$$ We have therefore recovered \eqref{thd}, and thus reduced the computation to the case we already studied. 

\begin {theorem}\label{prop2} Fix $r$ and $d$ two integers with $h=\gcd (r, d)$, and $N$ a vector bundle of slope $$\mu(N)=-\frac{d}{r},$$ and of minimal rank.  Then, $\mathsf E^{N}_{r, k}$ splits as sum of line bundles iff $h$ divides $k$. If $k=h(q-1)$, then $$\mathsf E^{N}_{r, k}\cong\left(\Theta_{1,(\det N)^{h}}\right)^{q-1}\otimes \left(\bigoplus_{\xi\in \widehat{\mathsf X}_h} \mathsf L_{\xi}^{\oplus \mathsf m_{\xi}}\right).$$  Here $\mathsf m_{\xi}$ are given by the same formulas \eqref{f1} and \eqref{f2} as in Theorem \ref{prop1}. 

\end {theorem}

\proof When $N_{0}=\mathsf W_{r', -d'}$, the statement is a consequence of the above discussion and the proof of Theorem \ref{prop1}. The general case follows  from here, since both the Verlinde bundle and the right hand side only change by translations. To see this, set $$L=\det N\otimes \left(\det N_{0}\right)^{-1}.$$ On the one hand, formulas of Drezet-Narasimhan \cite {dn} imply that  
$$\mathsf E^{N}_{r, k}={\det}_{\star}\left(\Theta_{r, N}^{k}\right)={\det}_{\star}\left(\Theta_{r, N_{0}}\otimes {\det}^{\star} L\right)^{k}=\mathsf E^{N_{0}}_{r, k}\otimes L^{k}.$$ In the above, we view the degree $0$ line bundle $L$ on $X$, as a line bundle on the Jacobian in the standard way. On the other hand, we have $$\Theta_{1, (\det N)^{h}}=\Theta_{1, (\det N_{0})^{h}}\otimes L^{h}.$$  The Theorem follows by putting these observations together. 
\qed

\subsection{Arbitrary level and rank.}\label{four} In this subsection we will prove Theorem \ref{prop3}. We will keep the same notations as in the introduction, writing $hr$ for the rank, and letting $hk$ be the level, with $\gcd (r, k)=1$. We will determine the splitting type of the Verlinde bundle $$\mathsf E_{hr, hk}={\det}_{\star}\left(\Theta^{hk}_{hr}\right)=\text {Sym}^{hk}\, \mathsf W_{hr},$$ obtained by pushing forward tensor powers of the canonical Theta bundle $\Theta_{hr}$ via $$\det:U_{X} (hr, 0)\to \text {Jac} (X)\cong X.$$ The case of non-zero degree and arbitrary Theta's is entirely similar, and we will leave the details to the interested reader. \vskip.07in

{\it Proof of Theorem \ref{prop3}.} We first consider the case when $r$ is odd. The arguments used to prove Theorem \ref{prop1} go through with only minor changes. It suffices to check that the decomposition \eqref{dec}: $$\mathsf E_{hr, hk}\cong \bigoplus_{\xi} \mathsf W_{r, k, \xi}^{\mathsf m_{\xi}}$$ holds $\mathsf X_{hr}$-equivariantly, after pullback by the morphism $hr$. The pullback of the left hand side is evaluated $\mathsf G_{hr}$-equivariantly via \eqref{lat1}: \begin{equation}\label{onemore}(hr)^{\star}\mathsf E_{hr, hk}\cong(hr)^{\star}\text{Sym}^{hk} \mathsf W_{hr}\cong \mathcal O_{X}(hr[o])^{hk}\otimes \text{Sym}^{hk}\mathsf S_{hk}^{\vee}.\end{equation} For the right hand side, recall first that $\mathsf W_{r, k, \xi}$ has rank $r$ and determinant $\mathcal O_{X}(k[o])\otimes \xi$. By comparing ranks and degrees, we see that \begin{equation}\label{amb}\mathsf W_{r, k, \xi}\cong\mathsf W_{r,k}\otimes \mathsf L_{\chi}.\end{equation} Here $\mathsf L_{\chi}$ is any $hr$-torsion line bundle with $\mathsf L_{\chi}^{r}=\mathsf L_{\xi}.$ Note $\chi$ is uniquely defined only up to $r$-torsion line bundles. The ambiguity inherent in the choice of $\chi$ will be shown to be inessential later. Observe that the pullback $(hr)^{\star}\mathsf L_{\chi}$ is the trivial bundle, endowed with the $\mathsf X_{hr}$-character $\chi$. 

We will determine the pullback of $\mathsf W_{r,k}$ by the morphism $hr$. As a first step, we show that non-equivariantly \begin{equation}\label{pullatiyah}r^{\star}\mathsf W_{r,k}\cong \mathcal O_{X}(kr[o])^{\oplus r}.\end{equation} The ingredients needed for the proof of \eqref{pullatiyah} are found in Lemma $22$ of Atiyah's paper \cite {atiyah}. There, it is explained that all indecomposable factors of $r^{\star}\mathsf W_{r,k}$ have the same rank $r'$ and degree $k'$. Therefore, $$r^{\star}\mathsf W_{r,k}\cong \bigoplus_{i=1}^{r/r'}\mathsf W_{r', k'}\otimes M_{i}$$ for some line bundles $M_{i}$. In fact, examining Atiyah's arguments, one can prove a little bit more. Using \eqref{inv2}, we observe that $$r^{\star} \mathsf W_{r,k}\otimes r^{\star}\mathsf W_{r,k}^{\vee}\cong \bigoplus_{1}^{r^{2}}\mathcal O_{X}.$$ The above tensor product contains $\mathsf W_{r',k'}\otimes \mathsf W_{r',k'}^{\vee}\otimes M_{i}\otimes M_{j}^{-1}$ as a direct summand, for any $i$ and $j$.  Now, applying equation \eqref{inv2} again, we see that $$\mathsf W_{r',k'}\otimes \mathsf W_{r',k'}^{\vee}\cong \bigoplus_{\rho}L_{\rho},$$ the sum being taken over the $r'$-torsion points $\rho$. This clearly gives a contradiction, unless $r'=1$ and the bundles $M_{i}$ and $M_{j}$ coincide. In conclusion, we proved that \begin{equation}\label{pull1}r^{\star}\mathsf W_{r,k}\cong \oplus_{i=1}^{r}M,\end{equation} for a suitable line bundle $M$. Taking determinants we obtain that $$M\cong\mathcal O_{X}(kr[o])\otimes P,$$ for some $r$-torsion line bundle $P$. We claim that $P$ is symmetric, {\it i.e.} $(-1)^{\star}P\cong P$. When $r$ is odd, these two facts together imply that $P$ must be trivial, proving \eqref{pullatiyah}. The symmetry of $P$ is a consequence of \eqref{pull1} and of the symmetry of $\mathsf W_{r,k}$. Indeed, $$(-1)^{\star}\mathsf W_{r,k}\cong \mathsf W_{r,k},$$ as both bundles are simple, and have the same rank and determinant. 

Having established \eqref{pullatiyah}, we compute 
\begin{equation}\label{pull2}(hr)^{\star} \mathsf W_{r, k}\cong \mathcal O_{X}(hr[o])^{hk}\otimes R,\end{equation} where $R$ is an $r$-dimensional vector space. In fact, $R$ carries a representation of the Theta group $\mathsf G_{hr}$, such that the center acts with weight $-hk$. However, this does not determine the representation $R$ uniquely, not even as a representation of $\mathsf H_{hr}$. In fact, one can show that there are precisely $h^{2}$ representations $R_{i,j}$ of $\mathsf H_{hr}$ with central weight $-hk$ \cite {S}; they will be indexed by integers $i, j\in \zh \times \zh$.

To determine $R$, we will use the following commutative diagram
\begin{center}
$\xymatrix{0\ar[r] & \mathbb C^{\star}\ar[r] & \mathsf G_{hr}\ar[r] &\mathsf X_{hr}\ar[r]&0\\ 0\ar[r] & \mathbb C^{\star}\ar[r]\ar[u]^{i} & \mathsf G_{h}\ar[u]\ar[r]& \mathsf X_{h}\ar[r]\ar@{_{(}->}[u] & 0.
}$
\end{center} Here, the morphism $i$ is the $r$-fold cover $\alpha\mapsto \alpha^{r},$ and the middle arrow is the natural morphism of Theta groups $\mathsf G_{h}\to \mathsf G_{hr}$. Via this diagram, we may consider the action of the group $\mathsf G_{h}$ on both sides of \eqref{pull2}. Recall from equation \eqref{dets} that $\mathsf G_{h}$-equivariantly, we have $$\mathcal O_{X}(hr[o])^{hk}\cong\mathcal O_{X}(h[o])^{hkr}\cong h^{\star}\mathcal O_{X}([o])^{kr}\otimes \left(\Lambda^{h}\mathsf S_{h}\right)^{kr}.$$ Therefore, using \eqref{pull2},  we see that $\mathsf G_{h}$-equivariantly, $$R\otimes \left(\Lambda^{h}\mathsf S_{h}\right)^{kr}=h^{\star}\left(r^{\star}\mathsf W_{r,k}\otimes \mathcal O_{X}(-kr[o])\right).$$ Note that the left hand side is an $\mathsf X_{h}$-module, since the center of $\mathsf G_{h}$ acts trivially; to this end, recall that the morphism $i$ is an $r$-fold covering of the centers.  By equation \eqref{pullatiyah}, the right hand side is the pullback of a trivial vector bundle, carrying a trivial $\mathsf X_{h}$-action. Consequently, the $\mathsf X_{h}$-representation $R\otimes \left (\Lambda^{h}\mathsf S_{h}\right)^{kr}$ is trivial. 

This latter observation pins down the $\mathsf H_{hr}$-representation $R$. Let us again pick theta structures, identifying the Theta group $\mathsf H_{hr}$ with the Heisenberg. The characters of the $h^{2}$ representations $R_{i,j}$ were computed in Theorem $3$ in \cite {S}. There it was proved that the trace of $\eta=(\alpha, x, y)\in \mathsf H_{hr}\cong \mu_{2hr}\times \mathbb Z/hr\mathbb Z\times \mathbb Z/hr\mathbb Z$ equals 
\begin {equation} \text {Trace}_{R_{i,j}}({\eta})=\begin {cases} r \alpha^{-hk}\zeta^{ix+jy} & \text {if } (x, y)\in \mathsf X_{h}, \text {{\it i.e. }if } (x,y)\in r\mathbb Z/hr\mathbb Z\times r\mathbb Z/hr\mathbb Z,\\ 0 & \text {otherwise.}\end{cases} \end {equation} Here $\zeta=\exp\left(\frac{2\pi i}{hr}\right).$  The character of the $\mathsf H_{h}$-representation $\left(\Lambda^{h}\mathsf S_{h}\right)^{kr}$ was calculated in \eqref{signs}: $$\text {Trace } (\eta)=\alpha^{hkr} (-1)^{(h+1)(x+y)kr}.$$ Since $R\otimes \left (\Lambda^{h}\mathsf S_{h}\right)^{kr}$ is a trivial $\mathsf X_{h}$-module, we must have $i=j=\frac{hrk(h+1)}{2}$. Then, the trace of $R$ becomes \begin{equation}\label{trace} \text {Trace}_{R}({\eta})=\begin {cases} r \alpha^{-hk}(-1)^{(h+1)kr(x+y)} & \text {if } (x, y)\in \mathsf X_{h},\\ 0 & \text {otherwise.}\end{cases} \end {equation}

Making use of \eqref{onemore} and \eqref{pullatiyah}, we can now check that both sides of \eqref{dec} agree equivariantly after pullback by $hr$. It remains to prove that $\mathsf H_{hr}$-equivariantly: \begin{equation}\label{final}\text{Sym}^{hk} \mathsf S_{hr}^{\vee}\cong R\otimes \bigoplus_{\chi} \chi^{\oplus \mathsf m_{\chi}}.\end{equation} In this sum, the $\chi$'s are $h^{2}$ representatives of the characters of $\mathsf X_{hr}$, modulo those characters of $\mathsf X_{hr}$ which restrict trivially to the subgroup $\mathsf X_{h}\hookrightarrow \mathsf X_{hr}.$ Taking representatives is necessary to avoid repetitions. Indeed, by comparing characters, we see that $$R\otimes \chi\cong R$$ iff $\chi$ restricts trivially to the subgroup $\mathsf X_{h}$. This equation also takes care of the ambiguity seemingly present in the pullback of \eqref{amb} by $hr$.  Note moreover that each representative $\chi$ appearing in the sum \eqref{final} restricts to a well-defined character $\xi$ of $\mathsf X_{h}$, hence giving rise to an $h$-torsion line bundle $\mathsf L_{\xi}$ on $X$. We will write $\omega$ for the order of this line bundle. 

In \eqref{final}, the multiplicities $\mathsf m_{\chi}$ are claimed to have the expressions given in equation \eqref{mul1} of the Theorem. Checking \eqref{final} amounts to a character calculation. For the left hand side, the character was essentially computed in Lemma \ref{chr}. Going through the proof of the Lemma, we see that the trace of $\eta=(\alpha, x, y)\in \mathsf H_{hr}$ on $\text{Sym}^{hk} \mathsf S_{hr}^{\vee}$ is zero, unless $(x,y)$ is an $h$-torsion point, say of order $h/\delta$ in $\mathsf X_{h}$. In the latter case, 
\begin{equation}\label{trace1}\text {Trace } (\eta)=(-1)^{xyk(h+1)}\alpha^{-hk}\cdot \frac{r}{r+k}\binom{(r+k)\delta}{r\delta}.\end{equation}  It suffices to check that the formula $$\mathsf m_{\chi}=\frac{1}{2h^{3}r}\sum_{\eta=(\alpha, x,y)\in\mu_{2hr} \times \mathsf X_{h}}\text{Trace }_{\text{Sym}^{hk}\mathsf S_{hr}}(\eta)\cdot \text {Trace}_{R}(\eta)^{-1}\cdot {\chi}(\eta)^{-1}$$ yields the same answer as \eqref{mul1}. Substituting \eqref{trace} and \eqref{trace1}, and recalling that $\xi$ denotes the restriction of $\chi$ to $\mathsf X_{h}$, we obtain $$\mathsf m_{\chi}=\frac{1}{h^{2}}\sum_{\delta|h}\frac{(-1)^{(h+1)kr\delta}}{r+k}\binom{(r+k)\delta}{r\delta}\sum_{(x,y) \text { has order } h/\delta}\xi(x,y).$$ By Lemma \ref{nl}, this expression can be rewritten as $$\mathsf m_{\chi}=\sum_{\delta|h}\frac{(-1)^{(h+1)kr\delta}}{(r+k)\delta^{2}}\binom{(r+k)\delta}{r\delta}\left\{\frac{h/\omega}{h/\delta}\right\}.$$ This completes the proof when $r$ is odd. 

When $r$ is even, $k$ must be odd, since $\gcd (r, k)=1$. Therefore, the Theorem holds true for the Verlinde bundle $\mathsf E_{hk, hr}$. We will now use the level-rank symmetry of the Verlinde bundles under the Fourier-Mukai transform $$\mathsf E_{hr, hk}^{\vee}\cong \widehat {\mathsf E_{hk, hr}},$$ which was explained in item $\text {(iv)}$ of the introduction. We claim that the Atiyah bundles enjoy the analogous symmetry under Fourier-Mukai: $$\mathsf W_{r, k, \xi}^{\vee}\cong \widehat{\mathsf W_{k, r, \xi}}.$$ Indeed, the case of trivial $\xi$ is the following well-known isomorphism generalizing \eqref{fmt}: $$\mathsf W_{r, k}^{\vee}\cong \widehat{\mathsf W_{r,k}}.$$ This is a consequence of the fact that both bundles are simple, of the same rank, and same determinant; alternatively, one may argue using the construction of the Atiyah bundles as successive extensions,  explained in Section \ref{geom}. The case of general $\xi$ is an immediate corollary, since the bundles involved differ only by translations. To see this, pick any line bundle $M$ with $M^{k}=\xi,$ and let $\tau_{M}$ denote the translation induced by $M$ on the elliptic curve. We compute 
$$\widehat {\mathsf W_{k, r, \xi}}\cong\widehat {\mathsf W_{k, r}\otimes M}\cong\tau_{M}^{\star}\widehat {\mathsf W_{k,r}} \cong \tau_{M}^{\star} W_{r, k}^{\vee}\cong \mathsf W_{r, k, \xi}^{\vee}.$$ The first and last isomorphism follow as usual by Atiyah's classification, while the second is a general fact about the Fourier-Mukai transform \cite {mukai}. 

We conclude the proof of the Theorem by collecting the above observations. 
\qed

\end{document}